\documentstyle[11pt]{article}
\def\mysection#1{\refstepcounter{section}\subsection{#1}}

\def\thefigure{\thefigure.\arabic{equation}}

\def\ns{\normalsize}
\def\ssize{\normalsize}
\def\note#1{}
\def\proof{\goodbreak\noindent{\bf Proof\quad}}
\def\endproof{{\ $\lform$}\bigskip }
\def\lform{\hbox{$\sqcup$}\llap{\hbox{$\sqcap$}}}
\def\ds{\displaystyle}
\def\ts{\textstyle}

\def\und#1{{\underline{#1}}}

\def\o{{}_{(1)}}\def\t{{}_{(2)}}
\def\Bo{{}_{\und{(1)}}}\def\Bt{{}_{\und{(2)}}}
\newtheorem{defn}{Definition}[subsection]
\newtheorem{thm}[defn]{Theorem}
\newtheorem{lm}[defn]{Lemma}
\newtheorem{prop}[defn]{Proposition}
\newcommand{\nn}{\nonumber}

\def\tens{{\mathop \otimes}}
\def\al{\alpha}
\def\be{\beta}
\def\del{\delta}

\def\ga{\gamma}
\def\la{\lambda}

\def\rh{\rho}

\def\th{\theta}
\def\eps{\epsilon}
\def\veps{\varepsilon}

\def\Om{\Omega}
\def\l{\left}
\def\r{\right}
\def\Z{{Z\kern-.647em Z}}
\def\CR{{\cal R}}
\def\CD{{\cal D}}
\def\oti{\otimes}
\def\tri{\triangle}
\def\qau{U_q\left(H_1,H_2,X^\pm\right)}
\def\un{{1}}
\def\ki{K_i}
\def\kii{K_i^{-1}}
\def\xpm{X^\pm}
\def\xp{X^+}
\def\xm{X^-}
\def\ko{K_1}
\def\kt{K_2}
\def\koi{K_1^{-1}}
\def\kti{K_2^{-1}}
\def\ho{H_1}
\def\htw{H_2}
\def\hi{H_i}
\def\dif{q-q^{-1}}
\def\fact{e^{i{\pi\over 2}\htw}}
\def\facti{e^{-i{\pi\over 2}\htw}}

\def\CRi{\CR^{-1}}
\def\quar{{1\over 4}}
\def\super{U_qgl\l(1|1\r)}
\def\vecl{{\bf l}}
\def\vect{{\bf t}}
\def\vecm{{\bf m}}
\def\vecu{{\bf u}}
\def\etap{\eta^+}
\def\cocross{{>\!\!\!\triangleleft\,}}
\def\lact{{\triangleright}}\def\ra{{\triangleleft}}
\def\extd{{{\rm d}}}

\begin{document}

\begin{titlepage}
\rightline{DAMTP/91-47/Revised}
\rightline{NIKHEF 95-023}
\vskip 1.8  true cm
\begin{center}
\Large{\bf Non-standard Quantum Groups and Superization}
\\
\vspace{0.45in}
\ns\sc
Shahn Majid
\footnote{Royal Society University Research Fellow and Fellow of
Pembroke College, Cambridge}\\
\ssize\em
Department of Applied Mathematics and
Theoretical Physics\\
University of Cambridge, Silver Street\\
Cambridge CB3 9EW, UK\\
\vspace{0.2in}
\ns\sc
M.J. Rodr\'\i guez-Plaza
\\
\ssize\em
NIKHEF-H\\
Postbus 41882\\
1009 DB Amsterdam, The Netherlands\\
\end{center}
\vskip 0.8 true cm
{\leftskip=1.5 true cm \rightskip=1.5 true cm
\noindent
We obtain the universal R-matrix of the non-standard quantum group
associated to the Alexander-Conway knot polynomial. We show further
that this non-standard quantum group is related to the super-quantum
group $U_qgl(1|1)$ by a general process of superization, which we
describe. We also study a twisted variant of this non-standard quantum
group and obtain, as a result, a twisted version of $U_qgl(1|1)$ as a
$q$-supersymmetry of the exterior differential calculus of any quantum
plane of Hecke type, acting by mixing the bosonic $x_i$ co-ordinates
and the forms $\extd x_i$.
\par}
\vskip .8 true cm
\end{titlepage}
\setcounter{page}{2}
\openup2\jot

\mysection{Introduction}
\label{intro}
Standard quantum group deformations $U_qg$ have been studied
extensively in recent years as deformations of the enveloping algebras
of ordinary Lie algebras \cite{Dri} \cite{Jim:dif}. One has techniques
for their construction coming from quantum and classical inverse
scattering (such as the FRT description \cite{FRT:lie} and the quantum
double \cite{Dri} of Drinfeld which provides their universal R-matrix
or quasitriangular structure).  One also has the possibility of
geometrical applications in analogy with the role of the corresponding
Lie algebras in the undeformed case.

Quite mysterious however, remain the quantum groups which one can
build in association to other, non-standard, solutions $R$ of the
quantum Yang-Baxter equations (QYBE)
$R_{12}\,R_{13}\,R_{23}=R_{23}\,R_{13}\,R_{12}$. The solutions of
these equations form a variety which has an interesting structure
\cite{Ma:ope} containing much more than the standard solutions and
their twistings. In low dimensions its structure consists
\cite{Hie:all} of several families, and the corresponding matrix
bialgebras $A(R)$ in these low-dimensional cases are computed
\cite{Hla:unu} but remain little understood. For example, existence of
a dual-quasitriangular structure for general $A(R)$ is known from
\cite{Ma:qua}, hence existence, at least formally, of a universal
R-matrix for
the corresponding `enveloping algebra' dual to it is also known, but
there is no easy algorithm to compute it (in general there may be
nothing like a borel subalgebra and associated quantum double as there
is in the standard case). The physical significance of such
non-standard quantum groups is likewise not at all understood because
of the absence of the analogies with undeformed objects which are possible
in the standard case.

In this paper, we study the simplest and most well-known of these
non-standard quantum groups in some detail.  As a Hopf algebra, it was
studied in \cite{JGW:new} but its full structure, notably the
universal R-matrix of its enveloping-type algebra
$U_q(H_1,H_2,X^\pm)$, has not been found before the announcement of
the present paper in \cite{MaPla:quan}.  We begin in
Section~\ref{quasi} by giving the full derivation of this, which is by
directly solving the non-linear equations for the quasitriangular
structure. It is, to our knowledge, the first non-standard quantum
group with non-trivial universal R-matrix known explicitly.  Many of
the standard constructions for quantum groups depend critically on the
universal R-matrix and our result means that these now hold
automatically for this non-standard quantum group. These include the
tensor product of tensor operators for the quantum group
\cite{Ma:qua}, quantum traces, and link and 3-manifold invariants. The
link invariants here depend on the ribbon element, which we compute,
and traces in irreducible representations. The matrix $R$ underlying
this example is known to be connected with the Alexander-Conway knot
polynomial \cite{CFFS:som}, which is usually developed in connection
with free fermions and the super-quantum group $U_qgl(1|1)$
\cite{KauSal:fre} \cite{BBCH:sym}, as well as in connection with other
quantum algebras \cite{GomSie:har}.

In Section~\ref{connection}, we make fully precise the connection
suggested here between this non-standard quantum group and
$U_qgl(1|1)$, through a process which we call {\em superization}.
This is the second goal of
the paper. This depends on a general algebraic principle of
transmutation introduced (by the first author) in
\cite{Ma:tra} as a way to shift the category (in this case
bosonic or super) in which an algebraic structure lives. We show that
this procedure agrees, for certain R-matrices, with a more ad-hoc
process of superization in which \cite{KulSkl:sol} certain types of
solutions of the QYBE can be converted to solutions of the super QYBE
(SYBE). We see this for our example. This elucidation is important
because it tells us that certain families of non-standard quantum
groups, corresponding to certain curves in the Yang-Baxter variety of
solutions, correspond after systematic superization to q-deformations
of super-Lie algebras. Without superization, we would have to view
super-Lie algebras as some independent theory `analogous' to the
standard theory of $U_qg$, but we see now that this theory is precisely
equivalent to the theory of certain non-standard bosonic quantum
groups in the same framework as their more standard cousins. One can
take this further and look for other more novel algebraic objects
corresponding to still other families of non-standard quantum groups
by further transmutation procedures. A step in this latter direction
is the sequel to the present paper \cite{MaPla:any}.

In Section~\ref{determinant}, we show how transmutation works for
the corresponding quantum groups of function algebra type. Again, there
is a general theory and a matching of our algebraic constructions with
the more ad-hoc process of replacing certain R-matrices by super
versions. We demonstrate the construction on the non-standard quantum
function algebra \cite{JGW:new} corresponding to the Alexander-Conway
R-matrix as in the previous section. This time, superization gives the
super-quantum group $GL_q(1|1)$ as studied for example in
\cite{SSV:pro}. Thus, once again, superization systematically relates
two different streams in the literature.  In the present case we
obtain, as an example, the (non-central) $q$-determinant for the
non-standard quantum group, corresponding to the known super
$q$-determinant of $GL_q(1|1)$.  The same relation holds between
certain non-standard quantum groups and $GL_q(n|m)$, among other
examples.

In Section~\ref{supersymmetry}, we proceed to a new application of (a
twisting of) our particular algebra $U_q(H_1,H_2,X^\pm)$ and its
superization, namely to $q$-deformed geometry.  A great many papers
have been devoted in recent years to the study of the exterior algebra
$\Omega_q$ of non-commuting `co-ordinates' $x_i$ and differential
forms $\extd x_i$ \cite{PusWor:twi} \cite{WesZum:cov} on quantum
planes. We show that, at least in the case when the exchange relations
governing the quantum plane are of Hecke type (e.g. the standard
$su_n$ family of quantum planes, but others as well), the exterior
algebra of differential forms has as a novel $q$-symmetry a
non-standard quantum group $U_q^\Omega(H_1,H_2,X^\pm)$ obtained by
twisting $U_q(H_1,H_2,X^\pm)$ by a quantum cocycle in the sense of
\cite{Dri:qua} \cite{GurMa:bra}.  The systematic theory of
superization tells us, further, that this corresponds equivalently to
a $q$-supersymmetry $U^\Omega_qgl(1|1)$ obtained by a super-cocycle
twist of $U_qgl(1|1)$. This is a general result which works for all
Hecke $q$-spaces of any dimension, with the odd generators acting by
mixing the $x_i$ and the $\extd x_i$.  As such, it generalizes an
important classical supersymmetry of the exterior algebra of
differential forms which has been used (for example) in \cite{HPPT}
and \cite{Wit:mor}, where it figures in a fundamental way. This
$q$-deformed $gl(1|1)$ supersymmetry of $q$-differential calculi does
not appear to be discussed elsewhere, and is a useful outcome of our
superisation techniques. It may also be possible to obtain its matrix
superalgebra form as a quotient of a general `universal
superbialgebra' construction \cite{LyuSud:sup}, which has been pointed
out to us.

This paper is the final version of a 1991 preprint with similar
title. We have added Section~\ref{supersymmetry} with the new
application to $q$-geometry. The earlier sections appear to us to
remain of interest, as well as playing a role in
\cite{CouLev:gen} and other subsequent works.

\mysection{Quasitriangular Structure on $U_q(H_1,H_2,X^\pm)$}
\label{quasi}

This section studies the Hopf algebra of quantum enveloping
algebra type associated to the Alexander-Conway matrix solution
\begin{equation}
R=\pmatrix{q&0&0&0\cr
           0&1&q-q^{-1}&0\cr
           0&0&1&0\cr
           0&0&0&-q^{-1}}
\label{AC}
\end{equation}
of the QYBE. The main result is an explicit expression of its universal
$R$-matrix $\CR$ or quasitriangular structure obeying the axioms of
Drinfeld \cite{Dri}. From this we show further that the Hopf algebra
is ribbon as well. We also compute $\CR$ in all finite irreducible
representations of the algebra at generic $q$.

We begin by constructing the quantum `enveloping' algebra by analogy
with the FRT approach \cite{FRT:lie}. I.e. we consider two matrices of
generators $\vecl^\pm$ and quadratic relations
$R\,\vecl^\pm_2\,\vecl^\pm_1=\vecl^\pm_1\,\vecl^\pm_2\,R$,
$R\,\vecl^+_2\,\vecl^-_1=\vecl^-_1\,\vecl^+_2\,R$ in the standard
notation $\vecl^\pm_1=\vecl^\pm\oti{\rm id}$, $\,\vecl^\pm_2={\rm
id}\oti\vecl^\pm$. This forms a bialgebra $\widetilde{U(R)}$ with
$\tri\vecl^\pm=\vecl^\pm\oti\vecl^\pm$, $\veps(\vecl^\pm)={\rm id}$ as
usual. We now quotient by adding further relations, or equivalently,
by making an ansatz for the specific form of $\vecl^\pm$. There is no
algorithm for this (apart from some general remarks in
\cite{Ma:qua}) but the lower and upper triangular ansatz
\begin{eqnarray}
\vecl^+=\l(\begin{array}{cc}
\ko                     &    0               \\
\l(\dif\r)\xp           &    \kti
\end{array}\r),
\qquad\qquad
\vecl^-=\l(\begin{array}{cc}
\koi               &      -\l(\dif\r)\xm     \\
0                       &   \kt
\end{array}\r)
\label{ansatz}
\end{eqnarray}
works and gives as quotient of $\widetilde{U(R)}$ the Hopf algebra
$U_q(K_1,K_2,X^\pm)$ say, studied in \cite{JGW:new} as generated by
$\un,~\ki,~\kii,~\xpm$, $i=1,2$ and relations
\begin{eqnarray}
&&\ki\cdot\kii=\kii\cdot\ki=1,\qquad \l[\ko,\kt\r]=0,\nn\\
&&\ko\xpm=q^{\pm 1}\xpm\ko,\qquad \kt\xpm=-q^{\mp 1}\xpm\kt,
\label{alg}\\
&&\l[\xp,\xm\r]={\ko\kt-\koi\kti\over \dif},\qquad\l(\xpm\r)^2=0.\nn
\end{eqnarray}
Here $q$ is an arbitrary parameter. The coalgebra structure is given by
\begin{eqnarray}
&&\tri \ki=\ki\oti\ki,\qquad \tri \kii=\kii\oti\kii,\nn\\
&&\tri \xp=\xp\oti\ko+\kti\oti\xp,\qquad \tri \xm=\xm\oti\kt+\koi\oti\xm,
\label{coal}\\
&&\veps\l(\ki\r)=\veps\l(\kii\r)=1,\qquad\veps\l(\xpm\r)\nn=0
\end{eqnarray}
and the antipode $S$ by
\begin{equation}
\begin{array}{c}
S\l(\ki\r)=\kii,\qquad S\l(\kii\r)=\ki,\\[12pt]
S\l(\xp\r)=-q\,\koi\kt\xp,\qquad S\l(\xm\r)=q\,\ko\kti\xm.
\end{array}
\label{ant}
\end{equation}

This is essentially the Hopf algebra which we study in this section.
Notice that the relations $\l(\xpm\r)^2=0$ are highly suggestive of a
superalgebra with $\xpm$ odd and $\ko,\kt$ even elements. Yet this is
an ordinary Hopf algebra and not at all a super-quantum group because
to extend the algebra homomorphism $\tri$ consistently to products of
generators is necessary to use the bosonic manipulation $(a\oti
b)(c\oti d)=a\,c\oti b\,d$ for the tensor product of two copies. In
the super-quantum group case we would need $(a\oti b)(c\oti
d)=(-1)^{{\rm deg}(b)~{\rm deg}(c)} a\,c\oti b\,d$, which turns out to
be inconsistent with the relations (\ref{alg}). In other words, this
non-standard quantum group reminds us of a super-quantum group but is
an ordinary bosonic one. This is a puzzle that we address in
Section~\ref{connection}.

We now introduce new generators for the above non-standard quantum
`enveloping algebra', namely
\begin{equation}
\ko=q^{\ho/2},\qquad \kt=\fact\,q^{\htw/2}.
\label{def}
\end{equation}

\begin{defn} We denote by $\qau$ the above non-standard Hopf
algebra with these new generators $\{\ho,\htw,\xpm\}$ and the
corresponding relations
\begin{equation}
\begin{array}{c}\l[\ho,\htw\r]=0,\qquad\l[\ho,\xpm\r]=\pm 2\xpm,
\qquad \l[\htw,\xpm\r]=\mp 2\xpm,\\[12pt]
\l[\xp,\xm\r]=\ds{\ko\kt -\koi\kti\over\dif}, \qquad\l(\xpm\r)^2=0,
\end{array}
\label{newalg}
\end{equation}
The coalgebra and antipode for the new generators is
$\tri \hi=\hi\oti\un +\un\oti\hi$ and $\veps\l(\hi\r)=0$ $(i=1,2)$
$S\l(\hi\r)=-\hi$  (unchanged for the rest).
\end{defn}

This change to new primitive generators $H_i$ is familiar for the
standard quantum groups $U_qg$ and can be made precise using formal
power series in the deformation parameter \cite{Dri}. The novel feature
in our non-standard case, however, is the additional factor
$e^{i\pi\htw/2}$ in (\ref{def}) which requires more work to make
precise.
We proceed formally but note that such exponentials do have a
well-defined meaning for the operator representations of $H_2$ which
we consider.  An alternative is to  adjoin $g=e^{i\pi\htw/2}$ as a seperate
group-like element, with corresponding commutation relations.

We are now ready to obtain the quasitriangular structure for
this new Hopf algebra. Recall that a Hopf algebra $U$ is called {\em
quasitriangular} if there is an invertible element $\CR$ in
$U\oti U$ that obeys the axioms \cite{Dri}
\begin{equation}
\tri '\l(a\r)=\CR\tri\l(a\r)\CRi
\label{dri1}
\end{equation}
for all $a$ in $U$ and
\begin{equation}
\l(\tri\oti{\rm id}\r)\CR=\CR_{13}\CR_{23},\qquad
\l({\rm id}\oti\tri\r)\CR=\CR_{13}\CR_{12}.
\label{dri2}
\end{equation}
Here $\tri '$ is defined as $\tau\circ\tri$ where
$\tau(x\oti y)\mapsto y\oti x$. The comultiplication $\tri '$
gives a second Hopf algebra structure on $U$ in addition to $\tri$
and the `{\sl universal $R$-matrix}' $\CR$ intertwines
them. Equations (\ref{dri2}) are evaluated in $U^{\oti^3}$ and
when $\CR$ is expressed as the formal sum $\CR=\sum_i a_i\oti b_i$,
$\CR_{12}, \CR_{13}$ and $\CR_{23}$ are then given by
the standard notation $\CR_{12}=\sum_i a_i\oti b_i\oti\un$,
$\CR_{13}=\sum_i a_i\oti \un\oti b_i$ and
$\CR_{23}=\sum_i\un\oti a_i\oti b_i$.
Then $\CR$ satisfies quantum Yang-Baxter equation
\begin{equation}
\CR_{12}\CR_{13}\CR_{23}=\CR_{23}\CR_{13}\CR_{12}
\label{QYBE}
\end{equation}
in $U^{\tens 3}$ as one sees easily from (\ref{dri1}) and (\ref{dri2}).

\begin{thm}
The Hopf algebra $\qau$ is quasitriangular with universal $R$-matrix
\begin{equation}
\CR=e^{\ds{-i{\pi\over 4}\htw\oti\htw}}\,
q^{\,\ds{{1\over 4}\l(\ho\oti\ho-\htw\oti\htw\r)}}
\l(\un\oti\un+(1-q^2)\,E\oti F\r),
\label{universal}
\end{equation}
where $E$ and $F$ denote the elements $\kt\xp$ and $\kti\xm$
respectively. \end{thm}
\proof
We start proving Drinfeld's axioms (\ref{dri1}). For brevity we shall
use the shorthand notation $q^{H\oti H}$ to denote the combination
$e^{-i{\pi\over 4}\htw\oti\htw}\, q^{{1\over
4}\l(\ho\oti\ho-\htw\oti\htw\r)}$ in $\CR$ above.  It is clear that
(\ref{dri1}) is satisfied for each $a=\hi$ since $\tri
'\hi\,(=\tri\hi)$ commutes with $q^{H\oti H}$ and $E\oti F$.  It is
also satisfied for $a=E$ because
$\l(\ko\kt\oti E+E\oti\un\r)\,q^{H\oti H}=
q^{H\oti H}\,\l(\un\oti E+E\oti\koi\kti\r)$ and
$\l(\un\oti E+E\oti\koi\kti\r)\,\l(\un\oti\un+(1-q^2) E\oti F\r)
=\l(\un\oti\un+(1-q^2) E\oti F\r)\,\l(\un\oti E+E\oti\ko\kt\r)$,
which put together produce
$\tri 'E\cdot\CR=
\allowbreak\l(\ko\kt\oti E+
\allowbreak E\oti\un\r)
\allowbreak\,q^{H\oti H}
\allowbreak
\,\l(\un\oti\un+(1-q^2)E\oti F\r)=
\allowbreak
\CR\cdot\tri E$. In an analogous manner we see that the
relation is satisfied for $a=F$ as well so we
conclude that (\ref{dri1}) holds for the entire Hopf algebra.

We prove now the first half of
(\ref{dri2}). {}From $\l(\tri\oti {\rm id}\r)\,\hi\oti\hi=
\hi\oti\un\oti\hi+\un\oti\hi\oti\hi$ we have that
$\l(\tri\oti {\rm id}\r)\,q^{H\oti H}=q^{H\oti\un\oti H}\,
q^{\un\oti H\oti H}$, with $q^{H\oti\un\oti H}$ and
$q^{\un\oti H\oti H}$ here understood in obvious notation. In turn
\[
\l(\tri\oti {\rm id}\r)\,\l(\un\oti\un+(1-q^2)\,E\oti F\r)=
\l(\un^{\oti ^3}+(1-q^2)\, E\oti\ko\kt\oti F\r)
\l(\un^{\oti^3}+(1-q^2)\,\un\oti E\oti F\r)
\]
which gives
\[
\l(\tri\oti {\rm id}\r)\,\CR=q^{H\oti\un\oti H}
\allowbreak\l(\un^{\oti^3}+(1-q^2)\, E\oti\un\oti F\r)\,
\allowbreak q^{\un\oti H\oti H}
\allowbreak\l(\un^{\oti^3}+(1-q^2)\,\un\oti E\oti F\r)=
\CR_{13}\CR_{23},
\]
once the commutator
$q^{\un\oti H\oti H}\,\l(E\oti\ko\kt\oti F\r) =\l(E\oti\un\oti
F\r)\,q^{\un\oti H\oti H}$ is used. The remaining
equation in (\ref{dri2}) is proved with a similar procedure.  \endproof

This goes significantly beyond the scope of ref. \cite{JGW:new}, where
the quasitriangular structure $\CR$ is not investigated. It is crucial
for numerous applications. In particular, the quantities $u=\sum_i
S(b_i)\,a_i$, $S(u)$, $z=u\,S(u)=S(u)\,u$ and
$r=u\l(S(u)\r)^{-1}=\l(S(u)\r)^{-1}u$ are of special interest in the
general theory \cite{Dri:alm}.  Here the invertible element $u$
implements the square of the antipode in the form $uau^{-1}=S^2(a)$
for all elements $a$ of the quantum group, $z$ is necessarily central
because it commutes with any $a$, and $r$ is group-like, i.e.$\tri
r=r\oti r$, and satisfies $rar^{-1}=S^4(a)$. We have,
\begin{lm}
The elements $u,\, S(u),\, z,\, r$ in $\qau$ are given by
\begin{eqnarray}
&&u=e^{i{\pi\over 4}\htw^2}q^{-\quar\l(\ho^2-\htw^2\r)}
\l(\un+(1-q^2)\,\koi\kti FE\r),
\label{u}\\
&&S(u)=e^{i{\pi\over 4}\htw^2}q^{-\quar\l(\ho^2-\htw^2\r)}
\l(\un+(1-q^2)\,EF\ko\kt\r),
\label{Su}\\
&&z=e^{i{\pi\over 2}\htw^2}q^{-{1\over 2}\l(\ho^2-\htw^2\r)}
\l(\un+(1-q^2)\,\l(\ko\kt EF+\koi\kti FE\r)\r),
\label{z}\\
&&r=e^{-i\pi\htw}q^{-\l(\ho+\htw\r)}.
\label{r}
\end{eqnarray}
\end{lm}
\proof
First we evaluate the element $u$. The universal $R$-matrix
(\ref{universal}) can be written in a better way for the present
calculations as the serie
$\sum_{l,n=0}^\infty{(-1)^n\over l!n!}\l(\quar\ln q\r)^l
\allowbreak\l(i{\pi\over 4}+\quar\ln q\r)^n
\allowbreak\l(\ho^l\htw^n\oti\ho^l\htw^n\r)
\allowbreak\l(\un^{\oti^2}+(1-q^2) E\oti F\r)$. To the sum
$\sum_i S(b_i)\,a_i$ contribute the terms
$S\l(\ho^l\htw^n\r)\ho^l\htw^n$
and $S\l(\ho^l\htw^n F\r)\ho^l\htw^n E$, equal to
$(-1)^{l+n}\ho^{2l}\htw^{2n}$ and
$(-1)^{l+n+1}\ko\kt\l(\ho+2\r)^{2l}\l(\htw-2\r)^{2n}FE$,
respectively. Restoring from these contributions the unexpanded
form of $u$ what we obtain is the expression (\ref{u}).
Relation (\ref{Su}) is the result of acting with the antipode map on
(\ref{u}) and taking into account that $S$ is an antialgebra
homorphism so that $S(ab)=S(b)\,S(a)$.
The element $z$ in (\ref{z}) is the product of results (\ref{u})
and (\ref{Su}). Concerning the derivation of (\ref{r}) it requires the
inverse of $S(u)$ given by
\[\l(S(u)\r)^{-1}=e^{-i{\pi\over 4}\htw^2}
q^{\quar\l(\ho^2-\htw^2\r)}
\l(\un-(1-q^2)\,EF\koi\kti\r)\]
as one can check easily. The final expression in (\ref{r}) is
calculated then from (\ref{u}) and this previous result.
We complete in this way the proof of the lemma.
\endproof

The elements $\CR$ and $u$, etc., of a quasitriangular Hopf algebra
are the key to numerous applications of quantum groups, among them the
construction of link and 3-manifold invariants. In this context, a
quasitriangular Hopf algebra is called {\em ribbon} if there is a
central element $\nu$ in it such that $\nu$ satisfies
\cite{ResTur:rib}
\[\nu^2=u\,S(u), \quad \tri \nu=\l(\CR_{21}\CR\r)^{-1}\l(\nu\oti \nu\r),
\quad
\veps(\nu)=1, \quad S(\nu)=\nu,\]
where $\CR_{21}$ denotes $\tau\circ\CR=\sum_i b_i\oti
a_i$.

\begin{prop}
$\qau$ is a ribbon Hopf algebra. The element $\nu$ and its inverse are
given by
\begin{eqnarray}
&&\nu=e^{i{\pi\over 4}\htw^2}q^{-\quar\l(\ho^2-\htw^2\r)}\ko\kt
\l(\un+(1-q^2)\koi\kti FE\r)
\label{v}\\
&&\nu^{-1}=e^{-i{\pi\over
 4}\htw^2}q^{\quar\l(\ho^2-\htw^2\r)}\koi\kti
\l(\un-(1-q^2)\,\ko\kt FE\r).\nn
\end{eqnarray}
\end{prop}
\proof
It is easy to see that indeed $\nu^2=z$ given in (\ref{z}).
That $\tri \nu=\l(\CR_{21}\CR\r)^{-1}\l(\nu\oti \nu\r)$ can be checked with
(\ref{universal}) and (\ref{v}).
The remaining two properties derive from (\ref{coal}) and (\ref{ant})
in a straightforward manner. Finally, that $\nu\cdot \nu^{-1}=\nu^{-1}\cdot
\nu=\un$ is again a simple computation.
\endproof

About the general properties of $\qau$ we mention that other
elements of interest are the casimirs (central elements) given by
\[c_1^2=e^{\ds{i\pi\htw}}q^{\ds{\ho+\htw}},\qquad
           c_2=\xp\xm-{1\over 2}\l({\ko\kt-\koi\kti\over\dif}\r).\]
Notice that $r=c_1^{-2}$ so $r$ too is central in our case. This
indicates that the antipode has order 4, that is to say $S^4={\rm id}$.
The element $c_1\equiv e^{i{\pi\over 2}\htw}q^{\l(\ho+\htw\r)/2}$
is not central but anticommutes with $\xpm$.

We turn now to one of the original applications of the universal
$R$-matrix, namely to obtain matrix solutions of the QYBE.  Since
$\CR$ obeys (\ref{QYBE}) abstractly in the algebra, so does the matrix
of $\CR$ in any representation of the quasitriangular Hopf algebra.
In the R-matrix setting \cite{FRT:lie} we know that the bialgebra
$\widetilde{U(R)}$ is dually paired with the quantum matrix bialgebra
$A(R)$, which leads to the {\em canonical representation} $\rh$
defined by \cite{FRT:lie} \cite{Ma:qua}
\begin{equation}
\rh^i{}_j\,({l^+}^k{}_l)=R^i{}_j{}^k{}_l{},\qquad
\rh^i{}_j\,({l^-}^k{}_l)={R^{-1}}~^k{}_l{}^i{}_j{},
\label{canon}
\end{equation}
in terms of our R-matrix. In the case of the AC solution (\ref{AC})
this gives us a representation
\begin{equation}
\begin{array}{c}
\rh(\ko)=\l(\begin{array}{cc}
                   q&0\\
                   0&1
               \end{array}\r),\qquad
\rh(\kt)=\l(\begin{array}{cc}
                   1&0\\
                   0&-q
                \end{array}\r),\\[12pt]
\rh(\xp)=\l(\begin{array}{cc}
                   0&1\\
                   0&0\end{array}\r),
\qquad \rh(\xm)=\l(\begin{array}{cc}
                   0&0\\
                   1&0
                 \end{array}\r).
\end{array}
\label{canonical}
\end{equation}

Finite dimensional irreducible representations of the algebra
$U_q(K_1,K_2,X^\pm)$ have been studied in \cite{JGW:new}, where it is
shown that the only `faithful' ones among them are two-dimensional and
equivalent to the following
\begin{equation}
\begin{array}{c}
\pi(\ko)=\l(\begin{array}{cc}
                   \la_1&0\\
                   0&q^{-1}\la_1
               \end{array}\r),\qquad
\pi(\kt)=\l(\begin{array}{cc}
                   \la_2&0\\
                   0&-q\la_2
                \end{array}\r),\\ \label{irrep}  \\
\pi(\xp)=\l(\begin{array}{cc}
                   0&\l(\la_1\la_2-\la_1^{-1}\la_2^{-1}\r)/\l(\dif\r)\\
                   0&0\end{array}\r),
\qquad \pi(\xm)=\l(\begin{array}{cc}
                   0&0\\
                   1&0
                 \end{array}\r),
\end{array}
\end{equation}
where $\la_1,\la_2$ are any non zero complex constants such that
$\l(\la_1\la_2\r)^2\neq 1$. Here a representation is said `faithful'
in the spirit of \cite{JGW:new} if $\pi(a)=\pi(b)$ implies $a=b$ for
$a,b$ in the vector space spanned by $\{K_i,K_i^{-1},X^\pm\}$. Each
representation (\ref{irrep}) is labelled by two numbers
$\l(\la_1,\la_2\r)$ and the tensor product of two of them decomposes
in irreducible representations as the following law indicates
$$\l(\la_1,\la_2\r)\oti\l(\mu_1,\mu_2\r)=
\l(\la_1\mu_1,\la_2\mu_2\r)\oplus
\l(q^{-1}\la_1\mu_1,-q\la_2\mu_2\r).$$

We use these irreducible representations for $\qau$ as well,
parametrizing them as
\[\la_1=q^{m_1},\qquad\la_2=e^{i\pi m_2}q^{m_2},\]
where $m_1,\,m_2$ are such that
$q^{2(m_1+m_2)} e^{2\,i\pi m_2}\ne 1$. The corresponding representation
$[m_1,m_2]$ is
\begin{equation}
\pi(\ho)=\l(\begin{array}{cc}
 2m_1  &   0  \\
 0     &  2\l(m_1-1\r)
\end{array}\r),\qquad
\pi(\htw)=\l(\begin{array}{cc}
2m_2   &   0   \\
0   &   2\l(m_2+1\r)
\end{array}\r),
\label{irrep'}
\end{equation}
together with $\pi(\xpm)$  as in (\ref{irrep}).

For these representations we can state the following proposition
\begin{prop}
In the canonical representation $\rho=[1,0]$ we recover
$\rho\oti\rho\l(\CR\r)=R$, the AC solution (\ref{AC}) of
Section~\ref{quasi}.
In the general representation $\pi=[m_1,m_2]$ of $\qau$ we have that
$\pi\oti\pi\l(\CR\r)$ is again (\ref{AC}) with the substitution of $q$
by $e^{i\pi m_2}q^{m_1+m_2}$.
\end{prop}
\proof
This follows by direct computation.
\endproof

In the above proposition, the part concerning $\rho$ is to be expected
from the general theory of matrix quantum groups (see \cite{Ma:qua},
Sec. 3) which therefore serves as a check on our universal $R$-matrix
(\ref{universal}). With respect to $\pi\oti\pi\l(\CR\r)$ we see that
it does not provide new or unfamiliar solutions to the QYBE but a
reparametrization of the AC solution again, cf. remarks in
\cite{JGW:new} for other aspects of these representations.
On the other hand, there are certainly other
representations of $\qau$ (for example the infinite dimensional
left-regular one) on which our universal $\CR$ can provide new braid
group actions and corresponding invariants.

\mysection{Connection of $\qau$ with super $\super$}
\label{connection}

In the last section we pointed out that the relations
$\l(\xpm\r)^2=0$ are indicative of some kind of super-Hopf
algebra structure. Yet $\qau$ is an ordinary Hopf algebra and
therefore not a super one at all. This puzzle was raised in
\cite{JGW:new} and we give now some insight into it by means of the
transmutation theory of \cite{Ma:bra} \cite{Ma:tra}, which enunciates
that under suitable circumstances is possible to transform Hopf
algebras into super-Hopf algebras and vice-versa.  As an application
of this superization construction we prove here that the superization
of (a quotient of) $\qau$ coincides with the super-quantum group
$\super$. We also study this procedure more generally and connect it
with a procedure of superizing the R-matrix itself under certain
conditions.

We start with an algebraic {\em superization} procedure for ordinary
Hopf algebras as follows. It is a special case of a theory in
\cite{Ma:tra}.

\begin{prop}\label{superise} If $H$ is a Hopf algebra containing a
group-like
element $g$ such that $g^2=\un$, there is a super-Hopf algebra $\und
H$, its superization, defined as the same algebra and counit as $H$,
and the comultiplication, antipode (if any) and quasitriangular
structure (if any) of $H$ modified to
\begin{equation}
\und \tri h=\sum h_{(1)}\,g^{\ds{-{\rm deg}\,(h_{(2)})}}\oti
h_{(2)},
\qquad
{\und S}\l(h\r)=g^{\ds{\,{\rm deg}\l(h\r)}}\, S\l(h\r)
\label{trans}
\end{equation}
and
\begin{equation}
\und \CR=\CR_g^{-1}\sum_i a_i\,g^{\ds{-{\rm deg}\,\l(b_i\r)}}\oti
b_i.\label{transR}
\end{equation}
\end{prop}
Here $h$ denotes an arbitrary element of $H$ with
comultiplication $\tri h=\sum h_{(1)}\oti h_{(2)}$ in the standard
notation and $\CR=\sum_i a_i\oti b_i$ denotes the universal $R$-matrix
of $H$. The element $\CR_g={\ts{1\over 2}} (\un\oti\un+\un\oti
g+g\oti\un-g\oti g)=\CR_g^{-1}$ is the nontrivial universal $R$-matrix
of the Hopf algebra $\Z'_2=\{\un, g\}$ with relations $g^2=\un$, $\tri
g=g\oti g$, $\veps\l(g\r)=1$ and $S\l(g\r)=g$. When $h$ is regarded as
an element of $\und H$ its grading is obtained with the action of $g$
on $h$ in the adjoint representation, that is by
$g\,h\,g^{-1}=(-1)^{\ts{{\rm deg}\l(h\r)}}\,h$ on homogeneous
elements.  This superization procedure represents a particular case of
a general transmutation theory in \cite{Ma:tra} for a Hopf algebra
containing or being mapped from a general quasitriangular Hopf algebra
in the role of $\Z_2'$. Other examples of transmutation include
anyonization \cite{Ma:any} \cite{MaPla:any} and complete transmutation
\cite{Ma:bra}.

In order to superize $\qau$ in this way, we first note that the role
of $g$ is played by $e^{{i\pi\htw}/2}$ introduced in
(\ref{def}). This element $g=e^{{i\pi\htw}/2}$ has the property that
$g^2$ is central and group-like. Hence it is natural to impose the
relation $g^2=\un$ in the abstract algebra and consider the quotient
$\qau/g^2-\un$.  Moreover, from Section~\ref{quasi} we know that
$\rh\,(g^2)=\un$ in the induced canonical representation
(\ref{irrep'}), so this further quotient is consistent with the
canonical representation and the pairing with $A(R)$.  The element
$g$ here
commutes with $\ho,\,\htw$ and anticommutes with $X^\pm$. We have
\begin{prop}
\label{propsuperization}
Let $g=e^{{i\pi\htw}/2}$. The superization of the Hopf algebra
$\qau/{g^2-\un}$ is the super-Hopf algebra $\super/e^{2\pi iN}-\un$.
Here $\super$ is defined by generators $h,N$ even and
$\eta$, $\etap$ odd and relations
\begin{equation}
\begin{array}{c}
[N,\eta]=-\eta ,\qquad [N,\etap]=\etap,\\[12pt]
\l\{\eta,\etap\r\}=\ds{{q^h-q^{-h}\over \dif}},\qquad
\eta^2=0,\qquad \l(\etap\r)^2=0,
\end{array}
\label{superalg}
\end{equation}
and $h$ central. The supercoalgebra is given by
\begin{equation}
\begin{array}{c}
\und \tri h=h\oti \un+\un\oti h,\qquad
\und \tri N=N\oti\un+\un\oti N,\\[12pt]
\und \tri \eta=\eta\oti q^{h-N}+q^{-N}\oti\eta,\qquad
\und \tri \etap=\etap\oti q^N+q^{-h+N}\oti\etap\\[12pt]
\und \veps(h)=\veps(N)=\veps(\eta)=\veps(\etap)=0.
\end{array}
\label{supercoal}
\end{equation}
This $\super$ has a super-quasitriangular structure given
by the expression
\begin{equation}
\und\CR=q^{-{\ds\l(h\oti N+N\oti h\r)}}
\l(\un\oti\un+(1-q^2)\,q^N\,\eta\oti q^{-N}\,\etap\r).
\label{superuniversal}
\end{equation}
\end{prop}
\proof
We apply the above superization construction to the quotient
$\qau/{g^2-\un}$ which is itself a Hopf algebra since $g^2$ is central
and group-like. The super-quantum group that results then
is easily recognizable as $\super$ given as in
(\ref{superalg})-(\ref{supercoal}) if we redefine
the generators as follows
\[
h=\l(\ho+\htw\r)/2,\qquad N=\htw/2,\qquad \eta=\xp,\qquad \etap=\xm g
\]
(so that $q^h=q^{(\ho+\htw)/2}$ and $q^N=q^{\htw/2}$).
Notice that a direct consequence of this definition is that $h$ is
central as stated. Notice also that the supercomultiplication
(\ref{supercoal}) is an algebra homomorphism consistent with relations
(\ref{superalg}) provided that we use super manipulation. To see it
let us compute as an example
$\und\tri\eta^2=\eta\,q^{-N}\oti q^{h-N}\,\eta-
q^{-N}\,\eta\oti\eta\,q^{h-N}=0$. The remaining cases are operated in
a similar manner.

The super-quasitriangular structure (\ref{superuniversal}) follows from
(\ref{universal}) and the transformation (\ref{transR}), but appears
with an overall factor
$c\equiv{\CR_g}e^{-i\pi N\oti N}q^{h\oti h}$.
This factor is central in $\super\oti\super$ and satisfies
$\l(\tri\oti{\rm id}\r)c=c_{13}\,c_{23},\:
\l({\rm id}\oti\tri\r)c=c_{13}\,c_{12}$, so it can be disregarded from
the final expression of $\und\CR$ without loss of generality.
\endproof

Once the super-universal $R$-matrix is known is possible to find the
super-ribbon Hopf algebra structure of $\super$ whose explicit
calculation we omit because it follows the same steps as in the
$\qau$ case. Instead we consider another procedure to obtain
super-quantum groups based on a super version of the FRT
construction. This involves solutions of the super Yang-Baxter
equation (SYBE).  It is said that an invertible matrix $\und R$ in
${\rm End}\l(V\otimes V\r)$ is a {\em super R-matrix} if it obeys the
{\em null-degree condition}
\begin{eqnarray}
\und R^a{}_c{}^b{}_d=0\quad {\rm when}
\quad p(a)+p(b)-p(c)-p(d)\ne 0\,\, ({\rm mod}\, 2),
\label{null}
\end{eqnarray}
and $\und R$ is a solution of the SYBE \cite{KulSkl:sol}
\begin{eqnarray}
(-1)^{\ds{p(e)\l(p(f)+p(c)\r)}}
\,\und R^b{}_e{}^a{}_f\,\und R^i{}_k{}^f{}_c\,\und R^k{}_j{}^e{}_d=
(-1)^{\ds{p(r)\l(p(s)+p(a)\r)}}
\,\und R^i{}_p{}^b{}_r\,\und R^p{}_j{}^a{}_s\,\und R^r{}_d{}^s{}_c,
\label{SYBE}
\end{eqnarray}
where sum over repeated indices is understood. By $V$ we denote a
$\Z_2$-graded vector space of dimension $d=n+m$ with linear basis
$\{e_i\},\, i=1,\ldots, d$ and the assumption that all vectors $e_i$
are homogeneous of degree $p(i)\equiv {\rm deg}(e_i)=0$ when
$i=1,\ldots,n$ and degree $p(i)=1$ when $i=n+1,\ldots,m$.  The
particularization $p(i)=0$ for all $i$ in (\ref{SYBE}) gives the
ungraded Yang-Baxter equation of the previous section and makes
(\ref{null}) an empty relation. The condition (\ref{null}) demanded in
the super case is equivalent to taking $\und R$ as an even matrix when
the degree of each matrix element is given by the rule $p(\und
R^a{}_c{}^b{}_d)=p(a)+p(b)-p(c)-p(d)$.  There is no loss of generality
in adopting this null degree assumption and, on the contrary, the
advantage that is reduces the calculations considerably.  Neither of
these two conditions is changed by a transformation in $V$ of the type
$e'_j=Q_j{}^i e_i$ with ${\rm deg}(e'_i)={\rm deg}(e_i)$ for all $i$,
i.e., both relations are invariant under a change of the form $\und
R\rightarrow \l(Q\oti Q\r)^{-1}\,\und R\l(Q\oti Q\r)$ with $Q$ any
non-singular $d\times d$ matrix of null degree.

It is clear that to any super $R$-matrix is possible to associate a
matrix super-bialgebra $\widetilde{U(\und R)}$ with generators $\un$
and $\l\{{m^\pm}{}^i{}_j\r\}\, i,j=1,\ldots,d$ and algebra and
coalgebra relations
\begin{equation}
(-1)^{\ds{p(e)(p(c)+p(f))}}
\,\und R^b{}_e{}^a{}_f\, {m^x}\,^f{}_{c}\,{m^y}\,^e{}_{d}
=(-1)^{\ds{p(r)(p(a)+p(s))}}
\,{m^y}\,^b{}_{r}\,{m^x}\,^a{}_{s}\,\und R^r{}_d{}^s{}_c,
\label{superFRTal}
\end{equation}
\begin{equation}
\und\tri {m^\pm}{}^i{}_{j}=
\sum_{k=1}^d {m^\pm}{}^i{}_{k}\oti {m^\pm}{}^k{}_{j},\quad
\und\veps\l({m^\pm}^i{}_j\r)=\delta^i{}_j,
\label{superFRTcoal}
\end{equation}
where $(x,y)=(+,+),\, (+,-),\, (-,-)$.  The generators of this
bialgebra are defined as homogeneous elements of degree
$p\l({m^\pm}\,^i{}_{j}\r)\equiv p(i)+p(j)\quad({\rm mod}\,2)$ and the
super-coalgebra structure (\ref{superFRTcoal}) satisfies the relation
(\ref{superFRTal}) provided that again to multiply two copies of the
algebra we use the super manipulation $\l(a\oti b\r)\l(c\oti
d\r)=(-1)^{{\ts{p(b)\,p(c)}}}\l(a\,c\oti b\,d\r)$.  These formulae
(\ref{superFRTal})-(\ref{superFRTcoal}) are the super version of the
usual formulae on \cite{FRT:lie}, as is clear on writing
(\ref{superFRTal}) in the compact form $\und
R\,\vecm^\pm_2\,\vecm^\pm_1=\vecm^\pm_1\,\vecm^\pm_2\,\und R$, $\und
R\,\vecm^+_2\,\vecm^-_1=\vecm^-_1\,\vecm^+_2\,\und R$. The notation is
similar to the ungraded version except that now the tensor product
contained in $\vecm^\pm_1, \,\vecm^\pm_2$ is $\Z_2$-graded. We
consider that the super version of the FRT method consists then in
finding suitable anstaze (or additional relations) for the $\vecm^\pm$
in order to obtain a super quantum group as quotient of
$\widetilde{U(\und R)}$. We consider now how this process is related
to the ungraded situation which we have considered before.

\begin{defn}\label{supableR} A solution $R\in M_d\tens M_d$ of
the usual QYBE is called {\em superizable} if there exists a grading
$p(i)\in\{0,1\}$ on the indices such that
\[ R^a{}_c{}^b{}_d=0\quad
{\rm when}\quad p(a)+p(b)-p(c)-p(d)\ne 0\ ({\rm mod} 2).
\]
\end{defn}

It is clear that superizable solutions of the QYBE correspond to super
R-matrices via the relation
\begin{equation}
\und R^a{}_c{}^b{}_d=(-1)^{\ds{p(a)p(b)}}\, R^a{}_c{}^b{}_d
\label{relation}
\end{equation}
We can build a super-bialgebra $\widetilde{U(\und R)}$ from this and ask
whether or not there is a reasonable ansatz for the super-matrix generators
$\vecm^\pm$  so that the resulting super-quantum group matches our algebraic
superization procedure starting from an ordinary quantum group built from
$\widetilde{U(R)}$.

\begin{prop} The AC R-matrix (\ref{AC}) is superizable with $p(1)=0,
p(2)=1$. The corresponding super R-matrix comes out from
 (\ref{relation}) as \begin{equation}
\und R=\pmatrix{q&0&0&0\cr
                0&1&q-q^{-1}&0\cr
                0&0&1&0\cr
                0&0&0&q^{-1}}.
\label{superAC}
\end{equation}
The super-bialgebra $\widetilde{U(\und R)}$ has as quotient the
quantum group $U_qgl(1|1)$ given by the ansatz
\begin{eqnarray*}
\vecm^+=\l(\begin{array}{cc}
q^{h-N}                     &    0               \\
\l(\dif\r)\eta           &    q^{-N}
\end{array}\r),
\qquad\qquad
\vecm^-=\l(\begin{array}{cc}
q^{-h+N}               &     -\l(\dif\r)\etap     \\
0                       &   q^N
\end{array}\r).
\end{eqnarray*}
\end{prop}
\proof We compute $\und R$ from (\ref{relation}) and then insert the
stated
ansatz into the relations (\ref{superFRTal})-(\ref{superFRTcoal}) and
recover those for $U_qgl(1|1)$ in Proposition~\ref{propsuperization}.
\endproof

Hence, at least in this case, our abstract superization construction as
in Proposition~\ref{superise} is compatible with other more ad-hoc
ideas of `superization' consisting of passing from R-matrices to super
R-matrices and looking for suitable ansatze on the quadratic FRT
bialgebras. It is possible to develop this second method further in
such a way as to always match with the algebraic superization. One has
to consider $p$ as defining a Hopf algebra homomorphism
$\widetilde{U(R)}\to \Z_2'$ as a dual-quasitriangular Hopf
algebra. This requires more machinery than we have introduced so far,
but is analogous to a slightly different consideration for matrix
quantum groups in \cite[Appendix]{Ma:varen}.

\mysection{Quantum determinant and antipode for $A(R)$ and its
connection with super $GL_q(1|1)$}
\label{determinant}

In this section we look at a non-standard quantum group built from
$A(R)$ where $R$ is our solution (\ref{AC}) and connect it by
superization with $GL_q(1|1)$ as studied for example in
\cite{SSV:pro}. This is dual to superization in the preceding section
but we will see that by adjoining an element with $g^2=1$ (rather than
quotienting by the relation $g^2=1$ as before) we can still apply our
algebraic superization theory (\ref{trans}).

Recall that $A(R)$ is the bialgebra generated by $\un$
and $\vect=\l\{t^i{}_j\r\}\, i,j=1,\ldots,d$ with algebra and
coalgebra relations given by $R\,\vect_1\,\vect_2=\vect_2\,\vect_1\,R$
and $\tri \vect=\vect\oti\vect$, $\veps\l(\vect\r)={\rm id}$,
respectively. In particular, for the AC solution (\ref{AC})  we have:

\begin{defn} $A(R)$ is the matrix bialgebra generated by
$\vect=\pmatrix{a&b\cr c&d}$ and relations
\begin{equation}
\begin{array}{c}
b\,a =q\, a\,b,\qquad c\,a = q\, a\,c,\qquad  d\,b =-q^{-1}\,b\,d,
\qquad d\,c=-q^{-1}\,c\,d,\\[12pt]
c\,b =b\,c,\qquad d\,a - a\,d=(q-q^{-1})\, b\,c, \qquad b^2= c^2=0,
\\[12pt]
\tri \vect=\vect\oti\vect,\qquad \veps(\vect)={\rm id}.
\end{array}
\label{A(R)}
\end{equation}
Assuming that $a$ and $d$ are invertible it is then possible to
define the antipode by
\begin{equation}
S(\vect)=\l(\begin{array}{cc}
a^{-1}+a^{-1} b d^{-1}c a^{-1}
&-a^{-1} b d^{-1}\\
 -d^{-1} c a^{-1}&
d^{-1}+d^{-1}c a^{-1} b d^{-1}
\end{array}\r),
\label{S(vect)}
\end{equation}
which makes the bialgebra $A(R)$ into a Hopf algebra.
\end{defn}
Unlike some other quantum groups, $SL(2)_q$ for example, the existence
of each generator antipode does not require any quantum determinant
condition. Furthermore, $S(\vect)$ can be computed much in analogy
with supermatrices as follows: the antipode map axioms demand that
$\vect\cdot S(\vect)=S(\vect)\cdot\vect=\un$ so in matrix terms we
must find a matrix $S(\vect)$ such that $S(\vect)=\vect^{-1}$. To
obtain it we first split $\vect$ in the sum of matrices
$\vect_0=\pmatrix{a&0\cr 0&d}$ and $\vect_1=\pmatrix{0&b\cr c&0}$ and
operate on
$\vect^{-1}=\l(\vect_0\l(\un+\vect_0^{-1}\vect_1\r)\r)^{-1}=
\l(\un-\vect_0^{-1}\vect_1+\l(\vect_0^{-1}\vect_1\r)^2\r)
\vect_0^{-1}$, where the power serie truncates because of relations
$b^2=c^2=0$. The only assumption for this to be
correct is that $a,\, d$ are invertible elements.

Regarding the existence of a quantum determinant in the Hopf algebra
there does not seem to be any group-like central element that could
play the role of quantum determinant properly. We explain this in more
detail noting that in $A(R)$ with $a,d$ invertible the element
\[\CD(\vect)=a d^{-1}-b d^{-1} c d^{-1}\]
commutes with $a,d$ and anticommutes with $b,c$, is invertible and
group-like. It also satisfies the relation $\CD(\vect\cdot
\vect')=\CD(\vect)\cdot\CD(\vect')$ for $\vect$ and $\vect'$ any two
commuting quantum matrices whose elements satisfy (\ref{A(R)}).  Of
course, this result means that $\CD^2$ is central and group-like so
that we could think of it as the quantum determinant of
$\vect$. However, we refrain from calling it the quantum determinant
since it is not particularly natural in this role.  For instance, the
relation $\CD^2=\un$ is not compatible with the duality pairing
$\l<\vecl^\pm,\vect\r>=R^\pm$ with $U_q(K_1,K_2,X^\pm)$ in the case of
solution (\ref{AC}) since the value of $\CD^2$ in the fundamental
$(+)$ and conjugate fundamental $(-)$ representation of $A(R)$ given
by
\begin{equation}
{\rh^+}^i{}_j\,(t^k{}_l)=R^k{}_l{}^i{}_j{},\qquad
{\rh^-}^i{}_j\,(t^k{}_l)={R^{-1}}~^i{}_j{}^k{}_l{}
\label{fund-conjfund}
\end{equation}
is $\rh^\pm\l(\CD^2\r)=q^{\pm 2}$. This indicates that
we could naturally set either $\CD^2=q^2$ or $\CD^2=q^{-2}$  not to
the unit.
The derivation of this relation from $R$ is
according to \cite{Ma:qua}, Sec 3.3.4.

These remarks about this non-standard quantum groups appear strange
from the point of view of quantum group theory but, once again, become
clear from the point of view of the corresponding super-quantum group
obtained by superization. To do this we apply the transmutation theory
of Section~\ref{connection} by first extending it by the group algebra
of $\Z_2$ as a Hopf algebra semidirect product $A(R)\cocross\Z_2$ and
then applying Proposition~\ref{superise}. This extension is nothing
but $A(R)$ with generators $a,b,c,d$ entirely unchanged and the extra
generator $g$ of $\Z_2$ adjoined, with $g^2=\un$, $\tri g=g\oti g$ and
the cross relations $a\,g=g\,a$, $b\,g=-g\,b$, $c\,g=-g\,c$,
$d\,g=g\,d$.  The product $A(R)\cocross\Z_2$ is still a Hopf algebra
because the comultiplications of $A(R)$ and $\Z_2$ extend as an
algebra homomorphism to $A(R)\cocross \Z_2$. In analogy with
Proposition~\ref{propsuperization} we have
\begin{prop}
\label{prosuperA(R)}
Let $A(R)$ be the matrix quantum group (\ref{A(R)}) associated to the
Alexander-Conway solution of the QYBE. The superization of its
$\Z_2$-extension $A(R)\cocross \Z_2$ is isomorphic to the super
$\Z_2$-extension $GL_q(1|1)\cocross \Z_2$.
Here $GL_q(1|1)$ is the super-Hopf algebra with generators
$\vecu=\pmatrix{\al&\be\cr \ga&\del}$, $\al$ and $\del$
even, $\be$ and $\ga$ odd and relations \cite{CFFS:som} \cite{SSV:pro}
\begin{eqnarray*}
\be\,\al=q\,\al\,\be,\qquad \ga\,\al=q\,\al\,\ga,
\qquad \del\,\be=q^{-1}\,\be\,\del,\qquad \del\,\ga=q^{-1}\,\ga\,\del,\\
\ga\,\be=-\be\,\ga,\qquad
\del\,\al-\al\,\del=(q-q^{-1})\,\be\,\ga\qquad \be^2=\ga^2=0,
\end{eqnarray*}
\[\und\tri\vecu=\vecu\oti\vecu,\qquad \und\veps(\vecu)={\rm id}\]
and
\begin{equation}
\und S(\vecu)=\l(\begin{array}{cc}
\al^{-1}+\al^{-1} \be\del^{-1}\ga\al^{-1}
&-\al^{-1} \be\del^{-1}\\
 -\del^{-1} \ga\al^{-1}&
\del^{-1}+\del^{-1} \ga\al^{-1} \be\del^{-1}
\end{array}\r).
\label{superS(vecu)}
\end{equation}
Its $\Z_2$-extension is by adjoining a bosonic element $g$
implementing the degree operator.
\end{prop}
\proof
{}From Section~\ref{connection} we know that the superization
of $A(R)\cocross \Z_2$ has the same algebra but a modified
comultiplication. Redefining the generators of this
Hopf algebra as
\begin{equation}
\al=a,\qquad \be=b\,g,\qquad \ga=c, \qquad \del=d\,g.
\label{redefinition}
\end{equation}
is immediate to see that in the new generators the algebra is that of
$GL_q(1|1)\cocross \Z_2$ while the superized comultiplication
becomes the matrix one of $GL_q(1|1)$ and $\und\tri g=g\oti g$
is unchanged. In $GL_q(1|1)\cocross \Z_2$ the
cross relations with $g$ are simply according to the super-statistics,
i.e. $g$ commutes with the even generators of $GL_q(1|1)$ and
anticommutes with the odd ones.
\endproof

We see in particular that the even combination
$\CD\,g=\al\del^{-1}-\be\del^{-1}\ga\del^{-1}$, group-like and central
in $A(R)\cocross \Z_2$, is viewed after superization as the usual
super-quantum determinant of $GL_q(1|1)$ as in
\cite{SSV:pro}.

We conclude this section with a general theorem of which the above
is an example. Its proof also supplies details of the proof of
Proposition~\ref{prosuperA(R)} previously sketched.
\begin{thm}
Let $R\in M_d\oti M_d$ be a superizable matrix solution of the
QYBE as in Definition~\ref{supableR} and let $A(R)$ be
the matrix FRT bialgebra associated to $R$.
Let $A(R)\cocross \Z_2$ be the $\Z_2$-extension of
$A(R)$ by adjoining an element $g$ with $g^2=\un$, $\tri g=g\oti g$
and the cross relations
\[
g\,t^i{}_j=(-1)^{\ds{p(i)+p(j)}}t^i{}_j\,g.
\]
Then $A(R)\cocross\Z_2$ is an ordinary bialgebra and its superization
is isomorphic to $A(\und R)\cocross \Z_2$, that is, to the
$\Z_2$-extension of the super FRT bialgebra associated to $\und R$.
\end{thm}
\proof
(i) First we check that we can make the $\Z_2$-extension of $A(R)$
as claimed. To do this we define an action of $\Z_2$ on $A(R)$ by
$g\,\lact t^i{}_j=(-1)^{\ts{p(i)+p(j)}}\,t^i{}_j$ extended to all
$A(R)$ as a module algebra, i.e.
$g\,\lact(t^i{}_j\,t^k{}_l)=(g\,\lact t^i{}_j)\,(g\,\lact t^k{}_l)$.
For this to be consistent with the algebra relations
$R\,\vect_1\,\vect_2=\vect_2\,\vect_1\,R$ of $A(R)$ we need
$(-1)^{\ts{p(f)+p(c)+p(e)+p(d)}}\, R^a{}_f{}^b{}_e\,t^f{}_c\,t^e{}_d=
(-1)^{\ts{p(b)+p(r)+p(a)+p(s)}}\,t^b{}_r\,t^a{}_s\, R^s{}_c{}^r{}_d$.
Since $R$ is superizable, the extra factors are both
$(-1)^{\ts{p(a)+p(b)+p(c)+p(d)}}$ so they cancel.
Hence the action of $g$ extends to $A(R)$ and respects its algebra
structure. It also respects the coalgebra structure because
$(g\oti g)\,\lact\tri t^i{}_j=
\sum_k(-1)^{\ts{p(i)+p(k)+p(k)+p(j)}}\,t^i{}_k\oti t^k{}_j=
\tri(g\,\lact t^i{}_j)$ as required. Since the action of $g$ thus
respects the algebra and coalgebra structure, it is immediate that
the semidirect product $A(R)\cocross \Z_2$ by this action is a
bialgebra.

\noindent
(ii) The super FRT bialgebra $A(\und R)$ with generators $\un$,
$\vecu=(u^i{}_j)$ of degree $p(u^i{}_j)=p(i)+p(j)$
is defined by relations
\begin{equation}
(-1)^{\ds{p(a)p(b)+p(c)p(e)}}
\,\und R^a{}_f{}^b{}_e\, u^f{}_{c}\,u^e{}_{d}
=(-1)^{\ds{p(c)p(d)+p(r)p(a)}}
\,u^b{}_{r}\,u^a{}_{s}\,\und R^s{}_c{}^r{}_d,
\label{superfunctalg}
\end{equation}
and $\und\tri\vecu=\vecu\oti\vecu$, $\veps(\vecu)={\rm id}$.
Its $\Z_2$-extension $A(\und R)\cocross\Z_2$ as a super-bialgebra
consists of adjoining a bosonic element $g$ with  $g^2=\un$,
$\und\tri g=g\oti g$ and cross relations
$g\,u^i{}_j=(-1)^{\ts{p(i)+p(j)}}u^i{}_j\,g$. We show that
$\th: A(R)\cocross \Z_2 \to A(\und R)\cocross \Z_2$
defined by $\th(t^i{}_j)=u^i{}_j\,g^{p(j)}$ and
$\th(g)=g$ is well defined as an algebra isomorphism (notice here
that transformation (\ref{redefinition}) is precisely of this type).
Indeed, applying $\th$ to both sides of
$R\,\vect_1\,\vect_2=\vect_2\,\vect_1\,R$ we find that
$
(-1)^{\ds{p(c)(p(e)+p(d))}}\,R^a{}_f{}^b{}_e
\,u^f{}_{c}\,u^e{}_{d}
\allowbreak
\,g^{\ts{p(c)+p(d)}}=
\allowbreak
(-1)^{\ds{p(r)(p(a)+p(s))}}
\allowbreak
\,u^b{}_{r}\,u^a{}_{s}
\,R^s{}_c{}^r{}_d
\allowbreak
\,g^{\ts{p(r)+p(s)}}$.
Since $R$ is superizable we can replace $g^{\ts{p(r)+p(s)}}$ on the
right by $g^{\ts{p(c)+p(d)}}$ and hence cancel it. Writing now $R$ in
terms of $\und R$ with (\ref{relation}) we obtain exactly the algebra
relations (\ref{superfunctalg}) of $A(\und R)$. The cross relations
also coincide, those in $A(\und R)\cocross\Z_2$ being given by
commutativity or anticommutativity according to the grading
$p(u^i{}_j)=p(i)+p(j)$. Hence $\th$ is an algebra map. The
superization theorem applied to $A(R)\cocross\Z_2$ does not change the
algebra structure, but it does change the comultiplication. The new
comultiplication from (\ref{trans}) is $\und\tri t^i{}_j=\sum_k
t^i{}_k\,g^{\ts{p(k)+p(j)}}\oti t^k{}_j$ and $\und\tri g=g\oti g$
(unchanged). {}From this it follows that $(\th\oti\th)\,\und\tri
(t^i{}_j\,g^{\ts{-p(j)}})=
\sum_k u^i{}_k\oti u^k{}_j=\und\tri\circ\th
(t^i{}_j\,g^{\ts{-p(j)}})$.  Thus, after superizing $A(R)\cocross
\Z_2$, the map $\th$ becomes an isomorphism of super bialgebras. This
completes the proof of the theorem.  Finally we remark that if
$A(R)\cocross \Z_2$ (or quantum group version of it) has an antipode
map $S$ then this is also superized and $\th$ induces an antipode
$\und S$ on the corresponding version of $A(\und R)\cocross\Z_2$
obeying $\und S(g)=g$ and $\sum_k\und S(u^i{}_k)\,u^k{}_j=\del^i{}_j=
\sum_k u^i{}_k\,\und S(u^k{}_j)$. This is obtained by applying $\th$.
\endproof

This theorem precisely connects the algebraic superization theory from
\cite{Ma:tra} as used in Section~\ref{connection} with the idea of
replacing $R$ by $\und R$ and making a `parrallel' super-FRT
construction. The algebraic method works more generally and applies
also to quotients of $GL$ versions of $A(R)$ provided the additional
relations are compatible with the $\Z_2$ action. We have given here
the most easily explained setting for superization in which the (cross
product) algebra does not change but the coalgebra is made super.  One
can also view the passage from $A(R)$ to $A(\und R)$ as a change of
product induced by a dual transmutation theory where $p$ provides a
map $A(R)\to\Z_2'$, as explained in \cite[Appendix]{Ma:varen}. From
this point of view $A(\und R)=B(R,Z)$ where $Z$ is a
super-transposition R-matrix defined by $p$ and $B(\ ,\ )$ is the more
general transmutation construction recently studied more (as quantum
braided groups) in \cite{Hla:bra} and elsewhere. See
\cite[Appendix]{Ma:varen} for details.

Let us note finally that there are many non-standard R-matrices beyond
the specific solution (\ref{AC}) on which we focused, and to which the
above superization constructions apply equally well. Thus the
enveloping algebra super-quantum group $U_qgl(n|m)$ and the matrix
super-quantum group $GL(n|m)_q$ are the superization of certain
non-standard quantum groups. These are associated by super-FRT type
constructions to the super R-matrices
\[
\und R_{gl(n|m)}=q\sum_{i=1}^n E_{ii}\oti E_{ii}
+q^{-1}\sum_{i=n+1}^{n+m}E_{ii}\oti E_{ii}
+\sum_{i\neq j} E_{ii}\oti E_{jj}
+(\dif)\sum_{j>i}(-1)^{\ds{p(i)p(j)}}\, E_{ij}\oti E_{ji}.
\]
where $E_{ij}$ here, for $1\le i,\, j\le n+m$, are the $(n+m)\times
(n+m)$-matrix given by $(E_{ij})^k{}_l=\del^k{}_i\del^l{}_j$ and
$p(i)$ is the function with values 0,1 depending on whether
$i=1,\ldots n$ or $i=n+1,\ldots n+m$, respectively. This super
R-matrix is the superization of the nonstandard solution of the
ordinary QYBE
\[R_{gl(n|m)}=q\sum_{i=1}^n E_{ii}\oti E_{ii}
-q^{-1}\sum_{i=n+1}^{n+m}E_{ii}\oti E_{ii} +\sum_{i\neq j}
(-1)^{\ds{p(i)p(j)}}\, E_{ii}\oti E_{jj} +(\dif)\sum_{j>i} E_{ij}\oti
E_{ji},\] of which (\ref{AC}) of previous sections constitutes the
particular case $n=m=1$. It leads by an FRT-type construction to
bosonic (not super) non-standard quantum groups which are, however,
strictly connected by superization with $U_qgl(n|m)$ and
$GL_q(n|m)$. In particular, $R_{gl(n|m)}$ do not directly generate
these super-quantum groups, as sometimes stated in the literature.

\mysection{$q$-supersymmetry of $q$-exterior algebras}
\label{supersymmetry}

In this section we develop a particular application of (a twisted
version of) the non-standard quantum enveloping algebra
$U_q(H_1,H_2,X^\pm)$ studied in Section~\ref{quasi}. We consider the
quantum exterior algebra $\Omega_q(R)$ on a quantum plane of any
dimension $n$, defined by generators $\l\{\pmatrix{x_1\cdots x_n\cr
\extd x_1\cdots \extd x_n}\r\}$ and relations
\begin{equation}
x_i\,x_j=q^{-1}\,x_b\,x_a R^a{}_i{}^b{}_j,\qquad
\extd x_i\, x_j=q\, x_b\,\extd x_a R^a{}_i{}^b{}_j,\qquad
\extd x_i\,\extd x_j=-q\, \extd x_b\,\extd x_a R^a{}_i{}^b{}_j,
\label{exterior}
\end{equation}
where $R$ is any $n^2\times n^2$  matrix solution of the Yang-Baxter
equations obeying
the Hecke condition $\l(P\,R-q)(P\,R+q^{-1}\r)=0$. Here $P$ denotes
the permutation matrix. The Hecke condition is sufficient to ensure
that relations (\ref{exterior}) are compatible with the usual graded
Leibniz rule and $\extd^2=0$. Exterior algebras of this (and more
complicated) form are quite well-studied now \cite{PusWor:twi}
\cite{WesZum:cov}. We have grouped the generators, however, as a
rectangular quantum matrix $A(R_\Omega:R)$ \cite{MaMar:glu}, which
are defined like the
usual FRT relations but with respect to two R-matrices. As explained
in \cite{Ma:varen}, this means automatically that there is a braided
addition law on $\Omega_q(R)$ (i.e., one does not need to show this, as
some authors have done) and, relevant to us, a coaction from the left
of the quantum group $A(R_\Omega)$. Here
\begin{equation}
R_\Omega=\pmatrix{q&0&0&0\cr
                     0&q&q-q^{-1}&0 \cr
                     0&0&q^{-1}&0\cr
                     0&0&0&-q^{-1}}
\label{ROmega}
\end{equation}
is a close relative of (\ref{AC}) and defines the non-standard matrix
quantum group $A(R_\Omega)$ with generators $1$,
$\vect=\pmatrix{a&b\cr c&d}$ and relations \cite{Ma:varen}
\begin{eqnarray*}
b\,a = a\,b,\qquad c\,a = q^2\,a\,c,\qquad  d\,b = -b\,d,
\qquad d\,c=-q^{-2}\, c\,d,\\
c\,b =q^2 b\,c,\qquad d\,a - a\,d=-(1-q^2)\,b\,c, \qquad b^2= c^2=0.
\end{eqnarray*}
Then $\Omega_q(R)$ is covariant under matrix multiplication from the
left by such quantum matrices (as well as being covariant from the
right under $A(R)$, which is the usual point of view).

This non-standard quantum group $A(R_\Omega)$ is different from, but a
close relative of, the non-standard quantum group studied in
Section~\ref{determinant}. Our first result makes this precise at
the dual level of the `enveloping' algebras.

\begin{prop} The QYBE solution $R_\Omega$ and the ansatz
\[
\vecl^+=\l(\begin{array}{cc}
q^{(\ho-\htw)/2}           &    0               \\
(\dif)\,q^{\ho/2}\,\xp           & \facti\,q^{(\ho-\htw)/2}
\end{array}\r),
\quad
\vecl^-=\l(\begin{array}{cc}
q^{-(\ho+\htw)/2}               &      -(\dif)\,q^{-\htw/2}\,\xm
\\
0                       &   \fact\,q^{(\ho+\htw)/2}
\end{array}\r)
\]
leads to the quantum enveloping algebra $U_q^\Omega(H_1,H_2,X^\pm)$
with the same algebra and counit as in (\ref{newalg})
in Section~\ref{quasi} but the coproduct and antipode
\[
\begin{array}{c}
\tri_\Omega \hi=\hi\oti\un+\un\oti\hi,\\[12pt]
\tri_\Omega \xp=\xp\oti q^{-\htw/2}+\facti\,q^{-\htw/2}\oti\xp,\quad
\tri_\Omega
\xm=\xm\oti\fact\,q^{(\ho+2\,\htw)/2}+q^{-\ho/2}\oti\xm,\\[12pt]
S_\Omega\l(\hi\r)=-\hi,\qquad
S_\Omega\l(\xp\r)=-q\,\fact\,q^{\htw}\,\xp,\qquad
S_\Omega\l(\xm\r)=q\,\facti\,q^{-\htw}\,\xm.
\end{array}
\]
This quantum group $U_q^\Omega(H_1,H_2,X^\pm)$ is a twisting in the
sense of \cite{Dri:qua} \cite{GurMa:bra} of the nonstandard quantum
group $U_q(H_1,H_2,X^\pm)$ of Section~\ref{quasi}, but with the
quantum cocycle
\[
\chi=q^{\ds{{1\over 4}\, H_2\oti (H_1+H_2)}}
\]
\end{prop}
\proof We chose the ansatz so that the algebra relations of
$\widetilde{U(R_\Omega)}$ recover the same relations (\ref{newalg})
as in Section~\ref{quasi}. The matrix coproduct
of $\vecl^\pm$ then determines the coproduct of $H_i,X^\pm$, which
comes out differently from Section~\ref{quasi}.  We can recognize it
as of the twisted form
$\tri_\Omega h=\chi(\tri h)\chi^{-1}$ for all $h$ in
$U_q(H_1,H_2,X^\pm)$. We check finally that $\chi$ here itself obeys
the quantum 2-cocycle condition $(1\tens\chi)({\rm id}\tens\tri)\chi=
(\chi\tens 1)(\tri\tens{\rm id})\chi$ and $(\eps\tens{\rm id})\chi=1$
as required for the twisting theory.
\endproof

So this non-standard quantum group, while not exactly the one in
Section~\ref{quasi}, is `gauge equivalent' to it in the sense of
twisting.  This means, for example, that we have at once its universal
R-matrix as $\CR_\Omega=\chi_{21}\,\CR\,\chi^{-1}$ in terms of the one
found in Section~\ref{quasi}, namely
\[
\CR_\Omega=e^{\ds{-i{\pi\over 4}\htw\oti\htw}}\,
q^{\,\ds{{1\over 4}\l((\ho-\htw)\oti(\ho+\htw)\r)}}
\l(\un\oti\un+(1-q^2)\,E\oti q^{-(\ho+\htw)/2}\,F\r)
\]
for the quantum group $U_q^\Omega(H_1,H_2,X^\pm)$.

\begin{prop} The non-standard quantum enveloping algebra
$U_q^\Omega(H_1,H_2,X^\pm)$
acts covariantly from the right on any
quantum exterior algebra $\Omega_q(R)$ associated to a Hecke solution
$R$ of the QYBE.  The action is
\[
x_i\ra\, \ho=2\,x_i,\quad x_i\ra\, \htw=0,
\quad x_i\ra\,\xp=q^{-1}\,\extd x_i,
\quad x_i\ra\,\xm=0,
\]
\[
\extd x_i\ra\, \ho=0,\quad
\extd x_i\ra\, \htw=2\,\extd x_i,\quad
\extd x_i\ra\,\xp=0,\quad \extd x_i\ra\,\xm=q\,x_i.
\]
\end{prop}
\proof The left coation of $A(R_\Omega)$ dualizes to a right action of
$U_q^\Omega(H_1,H_2,X^\pm)$.
We compute it by evaluating the matrix
coaction against the generators $\vecl^\pm$ in the usual way, which in
our case means in terms of the R-matrix $R_\Omega$. This then
gives the form stated for the action of our $H_1,H_2, X^\pm$
generators. Since the coaction is an algebra homomorphism it means, as one
can check explicitly, that the action is covariant in the sense
$(x\,y)\ra\, h=(x\ra\, h\o)(y\ra\, h\t)$ for all $x,y\in\Omega_q(R)$ and
$\tri_\Omega h=h\o\tens h\t$, say.
\endproof

So this non-standard quantum group has a geometrical role as a hidden
quantum group symmetry (valid even for $q=1$) of the exterior algebra
of a quantum plane, mixing $x_i$ and $\extd x_i$.

Finally, superization allows us to pass systematically to a super
version of these results.

\begin{prop} The superization of $A(R_\Omega)$ is the matrix
super-bialgebra
$A(\und R_\Om)$ with generators 1 and $\vecu=\pmatrix{\al&\be\cr
\ga&\del}$, $\al$ and $\del$ even, $\be$ and $\ga$ odd and relations
\begin{eqnarray*}
\be\,\al=\al\,\be,\qquad \ga\,\al=q^2\,\al\,\ga,
\qquad \del\,\be=\be\,\del,\qquad \del\,\ga=q^{-2}\,\ga\,\del,\\
\ga\,\be=-q^2\,\be\,\ga,\qquad
\del\,\al-\al\,\del=(1-q^2)\,\be\,\ga\qquad \be^2=\ga^2=0
\end{eqnarray*}
and matrix super coproduct. If we assume $\al$ and $\del$ invertible,
we obtain a super-quantum group $GL_q^\Omega(1|1)$, say, with antipode
given by the usual formula (\ref{superS(vecu)}). The super determinant
is again a central, bosonic and group-like element given by
\[\und\det(\vecu)=\al\del^{-1}-\be\del^{-1}\ga\del^{-1}.\]
\end{prop}
\proof The matrix $R_\Omega$ is superizable with $p(1)=0$, $p(2)=1$. We
compute $\und R_\Omega$ from (\ref{relation}) (only a sign changes)
and put it into (\ref{superfunctalg}). If we allow $a,d$ invertible we
can obtain a nonstandard quantum group with antipode, much as in
Proposition~\ref{prosuperA(R)}.  Superising this gives a superantipode
if we assume $\alpha,\delta$ invertible.
\endproof

Next, a general feature of the transmutation theory is that
representations and covariant systems under the original object
remain so (but in the new category) under the transmuted one. Hence in
particular, representations or covariant systems under the original
bosonic quantum group become automatically (in our case)
super-representations or super-covariant systems under the
superization. We check this explicitly for our example.

\begin{prop} Let $\Omega_q(R)$ be a quantum exterior algebra of Hecke type,
regarded as a super algebra with $x_i$ even and $\extd x_i$ odd. Then
(keeping in mind the corresponding bose-fermi statistics) it is
covariant under the $GL_q^\Omega(1|1)$ transformation
\[\label{supertrans} \pmatrix{x_1'\cdots x_n'\cr \extd x_1'\cdots \extd
x_n'}=\pmatrix{\alpha&\beta\cr\gamma&\delta}\pmatrix{x_1\cdots x_n\cr
\extd x_1\cdots \extd x_n}\] in the sense that the primed co-ordinates
and forms obey the same relations of $\Omega_q(R)$.
\end{prop}
\proof We verify directly that indeed $\extd x'_i\,\extd x'_j=-q\,
\extd x'_b\,
\extd x'_a\,R^a{}_i{}^b{}_j$. From the stated super-transformation we have
that $\extd x'_i\,\extd x'_j+q\, \extd x'_b\,\extd x'_a\,
R^a{}_i{}^b{}_j=
\l(\ga\, x_i+\del\,\extd x_i\r)\l(\ga\, x_j+\del\,\extd x_j\r)
+ q \l(\ga\, x_b+\del\,\extd x_b\r)\l(\ga\, x_a+\del\,\extd x_a\r)
R^a{}_i{}^b{}_j
=\ga\,\del \, x_i\extd x_j-\del\,\ga\,\extd x_i\, x_j+
q\,\ga\,\del \, x_k\extd x_l\, R^l{}_i{}^k{}_j
-q\,\del\,\ga\,\extd x_k\, x_l\, R^l{}_i{}^k{}_j
=\ga\,\del \, x_i\extd x_j-\del\,\ga \extd x_i\, x_j+
\ga\,\del\,\extd x_i\,x_j
-q\,\del\,\ga\, x_e\extd x_f R^f{}_k{}^e{}_l\, R^l{}_i{}^k{}_j=
0$. In the last equality we used the Hecke condition in the form
$R^f{}_k{}^e{}_l\, R^l{}_i{}^k{}_j=(q-q^{-1})\,R^f{}_i{}^e{}_j+
\del^f_j\,\del^e_i$. Similarly for the other relations of $\Omega_q(R)$.
\endproof

The supersymmetry in this form appears to be related to a somewhat
larger `universal superbialgebra' coacting on exterior algebras,
developed by other means in~\cite{LyuSud:sup}.

Finally, to complete our picture, we construct the corresponding
superenveloping algebra $U_q^\Omega gl(1|1)$ either by algebraic
superization as in Section~\ref{connection} or from $\und R_\Omega$
and the ansatz
\begin{eqnarray*}
\vecm^+=\l(\begin{array}{cc}
q^{h-2\,N}                     &    0               \\
\l(\dif\r)\,q^{h-N}\,\eta           &    q^{h-2\,N}
\end{array}\r),
\qquad\qquad
\vecm^-=\l(\begin{array}{cc}
q^{-h}               &     -\l(\dif\r)\,q^{-N}\,\etap     \\
0                       &   q^h
\end{array}\r).
\end{eqnarray*}

\begin{prop} The super-quantum group $U_q^\Omega gl(1|1)$ has generators
$h,\,N,\,\eta,\etap$ obeying the same algebra relations as $\super$
but different supercomultiplication, superantipode and super-universal
R-matrix
\begin{equation}
\begin{array}{c}
{\und\tri}_\Omega h=h\oti \un+\un\oti h,\qquad
{\und\tri}_\Omega N=N\oti\un+\un\oti N\\[12pt]
{\und\tri}_\Omega \eta=\eta\oti q^{-N}+q^{-N}\oti\eta,\qquad
{\und\tri}_\Omega \etap=\etap\oti q^{h+N}+q^{-h+N}\oti\etap
\end{array}
\end{equation}
\begin{equation}
\und S_\Omega(h)=-h,\qquad \und S_\Omega(N)=-N,\qquad
\und S_\Omega(\eta)=-q\,q^{2N}\,\eta,\qquad
\und S_\Omega(\etap)=-q\,q^{-2N}\,\etap,
\end{equation}
\begin{equation}
{\und{\CR}}_\Omega=q^{-{\ds 2\,N\oti h}}
\l(\un\oti\un+(1-q^2)\,q^N\,\eta\oti q^{-N-h}\,\etap\r).
\end{equation}
It is related by twisting of super-quantum groups to $U_qgl(1|1)$ via
$\und\chi=q^{N\oti h}$, viewed as a super quantum 2-cocycle. Moreover, it
acts covariantly on any $\Omega_q(R)$ of Hecke type by
\[ x_i\ra\, h=x_i,\quad x_i\ra\, N=0,
\quad x_i\ra\,\eta=q^{-1}\,\extd x_i,
\quad x_i\ra\,\eta^+=0\]
\[ \extd x_i\ra\, h=\extd x_i,\quad
\extd x_i\ra\, N=\extd x_i,\quad \extd
x_i\ra\,\eta=0,\quad \extd x_i\ra\,\eta^+=q\,x_i\]
\end{prop}
\proof The relations of $\widetilde{U(\und R_\Omega)}$ for
the ansatz shown leads to the usual relations of $U_qgl(1|1)$.  The
matrix supercoproduct of the $\vecm^\pm$ gives the form on the
generators. We recognize it as $\und \tri$ in (\ref{supercoal})
twisted as
${\und\tri}_\Omega(\cdot )=\und\chi\,\und\tri(\cdot )\,\und\chi^{-1}$
for $\und\chi$ as shown. It is as super 2-cocycle in the sense
 \[
(\un\oti\und\chi)\,({\rm id}\oti\und\tri)\,\und\chi=
(\und\chi\oti\un)\,(\und\tri\oti {\rm id})\,\und\chi
\]
and $(\und\veps\oti {\rm id})\,\und\chi=\un$. This then gives at
once the super universal $R$-matrix as
stated, obtained from (\ref{superuniversal}) as
${\und\CR}_\Omega=\und\chi_{21}\,\und\CR\,\und\chi^{-1}$.  Note that, in
the present case, $\und\chi$ is bosonic so it is an ordinary 2-cocycle
just as well; the bosonic twisting (as in Proposition~V.1) followed by
superisation gives the same answer as superising first and then twisting
as we do now.  Finally, we dualize the coaction in the preceding
proposition via the super duality relations
\begin{equation}
\begin{array}{c}
<u{}^i{}_j,\,{m^+}\,^k{}_l>=(-1)^{\ds{p(j)(p(k)+p(l))}}
\,\und {R_\Om}{}^i{}_j{}^k{}_l{},\\[12pt]
<u{}^i{}_j,\,{m^-}\,^k{}_l>=(-1)^{\ds{p(k)(p(i)+p(j))}}
\,\und {R_\Om}^{-1}{}^k{}_l{}^i{}_j{}
\label{duality}
\end{array}
\end{equation}
to obtain the right actions as stated.
It follows (as one can also verify directly) that this right action of
$U_q^\Omega gl(1|1)$ is super-covariant in the sense
\begin{equation} (x\,y)\ra\,h=(-1)^{\deg(y)\,\deg(h\Bo)}\,(x\ra\,
h\Bo)\,(y\ra\, h\Bt)
\end{equation}
for all $x,y\in \Omega_q(R)$ and $\und\tri_\Omega h=h\Bo\tens h\Bt$ the
super-coproduct. The action on products of generators of $\Omega_q(R)$ is
consistently determined by this super-covariance.
\endproof


%
\def\section{\subsection}

\end{document}